\documentclass[11pt]{article}
\usepackage{indentfirst}
\usepackage{tikz}
\usepackage{graphics}
\usepackage{graphicx}
\usepackage{epsfig}
\usepackage{epstopdf}
\usepackage{subfigure}
\usepackage{amssymb}
\usepackage{amsthm}
\usepackage{mathrsfs}
\usepackage{hhline}
\usepackage{longtable}
\usepackage{amsmath}
\usepackage{array}
\usepackage{float}
\usepackage{multirow}
\usepackage{cite}
\usepackage{enumerate}
\usepackage{xcolor}
\usepackage{caption}
\numberwithin{equation}{section}
\newtheoremstyle{de}
  {10pt}          
  {10pt}  
  {\rm}  
  {}
  {\bf}  
  {. }    
  { }    
  {}     
\theoremstyle{de}

\newtheorem{theorem}{Theorem}[section]

\newtheorem{remark}[theorem]{\bf{Remark}}

\newtheorem{example}{Example}[section]
\newtheorem{problem}{Problem}[section]
\newtheoremstyle{theorem}
  {10pt}          
  {10pt}  
  {\it}  
  {}
  {\bf}  
  {. }    
  { }    
  {}     
\theoremstyle{theorem}

\newtheorem{defi}{Definition}[section]

\usepackage{xcolor}

\definecolor{vividviolet}{rgb}{0.62, 0.0, 1.0}

\newcolumntype{"}{@{\hskip\tabcolsep\vrule width 1pt\hskip\tabcolsep}}
\numberwithin{equation}{section}
\def\ms{\medskip}
\def\nt{\noindent}

\graphicspath{%
    {converted_graphics/}
    {/}
}
\begin{document}
\baselineskip18truept
\normalsize
\begin{center}
{\mathversion{bold}\Large \bf An algorithmic approach in constructing infinitely many even size graphs with local antimagic chromatic number 3}

\bigskip
{\large Gee-Choon Lau}

\medskip

\emph{College of Computing, Informatics \& Mathematics,}\\
\emph{Universiti Teknologi MARA (Segamat Campus),}\\
\emph{85000, Johor, Malaysia.}\\
\emph{geeclau@yahoo.com}\\


\end{center}
%

\begin{center} Dedicated to Prof Arumugam S. on the occasion of his 80th birthday
\end{center}

\begin{abstract}
An edge labeling of a connected graph $G = (V, E)$ is said to be local antimagic if it is a bijection $f:E \to\{1,\ldots ,|E|\}$ such that for any pair of adjacent vertices $x$ and $y$, $f^+(x)\not= f^+(y)$, where the induced vertex label $f^+(x)= \sum f(e)$, with $e$ ranging over all the edges incident to $x$.  The local antimagic chromatic number of $G$, denoted by $\chi_{la}(G)$, is the minimum number of distinct induced vertex labels over all local antimagic labelings of $G$. In this paper, we first introduce an algorithmic approach to construct a family of infinitely many even size non-regular tripartite graphs with $t\ge 1$ component(s) in which every component, called a {\it Luv} graph, is of odd order $p\ge 9$ and size $q=n(p+1)$ for $n\ge 2$. We show that every graph in this family has local antimagic chromatic number 3. We then allowed the $m$-th component to have order $p_m\ge 9$ and size $n_m(p_m+1)$ for $n_m\ge 2, 1\le m\le t$. We also proved that every such graph with all components  having same order and size also has local antimagic chromatic number 3. Lastly, we constructed another family of infinitely many graphs  such that different components may have different order and size all of which having local antimagic chromatic number 3. Consequently, many other families of (possibly disconnected) graphs with local antimagic chromatic number 3 are also constructed.  

\ms\noindent Keywords: Local antimagic chromatic number, disconnected, non-regular, tripartite

\end{abstract}

\baselineskip18truept
\normalsize



\section{Introduction}

\nt Let $G=(V, E)$ be a connected graph of order $p$ and size $q$. A bijection $f: E\rightarrow \{1, 2, \cdots, q\}$ is called a \textit{local antimagic labeling}
if $f^{+}(u)\neq f^{+}(v)$ whenever $uv\in E$, where $f^{+}(u)=\sum_{e\in E(u)}f(e)$ and $E(u)$ is the set of edges incident to $u$.
The mapping $f^{+}$ which is also denoted by $f^+_G$ is called a \textit{vertex labeling of $G$ induced by $f$}, and the labels assigned to vertices are called \textit{induced colors} under $f$. The \textit{color number} of a local antimagic labeling $f$ is the number of distinct induced colors under $f$, denoted by $c(f)$.  Moreover, $f$ is called a \textit{local antimagic $c(f)$-coloring} and $G$ is local antimagic $c(f)$-colorable. The \textit{local antimagic chromatic number} $\chi_{la}(G)$ is defined to be the minimum number of colors taken over all colorings of $G$ induced by local antimatic labelings of $G$ (see~\cite{Arumugam, Bensmail}). Affirmative solutions on some problems raised in~\cite{Arumugam} can be found in~\cite{LNS}. Let $G+H$ and $mG$ denote the disjoint union of graphs $G$ and $H$, and $m$ copies of $G$, respectively. In~\cite{Haslegrave}, the author proved that every connected graph of order at least 3 has a local antimagic labeling. Thus, local antimagic chromatic number is well defined for all graphs without a $K_2$ component. For integers $a < b$, let $[a,b] = \{a, a+1, \ldots, b]$. For an edge $e$ of $G$, let $G-e$ be the graph with edge $e$ deleted. Very little is known about the local antimagic chromatic number of disconnected graphs (see~\cite{Baca, CLS, LS-AMH} for some results on 2-regular graphs and forests). 

\begin{figure}[H]
\begin{center}
\centerline{\epsfig{file=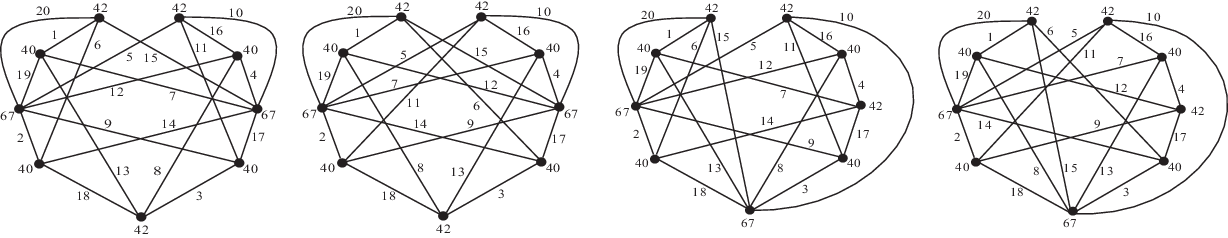, width=16cm}}
\caption{Four Luv graphs with local antimagic chromatic number 3.}\label{fig:4luv}
\end{center}
\end{figure}

\ms\nt Consider the following four non-isomorphic {\it Luv} graphs in Figure~\ref{fig:4luv}.  They are tripartite with exactly two degree 6 vertices and all others of degree 4. The associated edge-labeling shows that they have local antimagic chromatic number 3.

\ms\nt For $n\ge 2$, we can now define a family of non-regular disconnected tripartite graphs with $t\ge 1$ component(s) in which every component is a Luv graph of odd order $p\ge 9$ and size $q=n(p+1)$, having exactly $n$ degree $2n+2$ vertices and $p-n$ degree $2n$ vertices, as in Figure~\ref{fig:4luv}, denoted $\mathcal L_t(p,n)$ (thus known as  {\it Luv-All-in-Universe} graphs (abbreviated {\it LAU} graphs). The $m$-th component of each graph in $\mathcal L_t(p,n)$, denoted $L_t(p,n|m,a^1_{m},a^2_m,\ldots,a^n_m)$, has vertex set $V = \{u_{m,i}\,|\, 1\le m\le t, 1\le i\le p\}$ and edge set $E_m = E^\phi_{m} \cup F_m$ ($1\le \phi\le n$) such that 
\begin{enumerate}[(1)]
  \item $E^1_{m} = \{u_{m,i}u_{m,i+1}\,|\,1\le m\le t, 1\le i\le p-1\}$,
  \item $E^\phi_{m} = \{u_{m,j^\phi_i}u_{m,j^\phi_{i+1}}\,|\, 1\le i\le p-1, |j^\phi_{i+1} - j^\phi_i| \mbox{ is odd}\ge 3, 1=j^\phi_1 \ne j^\phi_2 \ne \cdots \ne j^\phi_{p-1} \ne j^\phi_p=p\}$ where $j^\phi_i, j^\phi_{i+1}$ are of distinct parity, and $E^\phi_m\cap E^{\phi\prime}_m = \varnothing$  for $2\le \phi \ne \phi\prime\le n$, 
  \item $F^\phi_m = \{u_{m,1}u_{m,2a^\phi_m+1}, u_{m,p}u_{m,2a^\phi_m+1}\}$ for $1\le a^1_m\ne \cdots \ne a^n_m\le \frac{p-3}{2}$. 
\end{enumerate} 
\nt Thus, the $p-1$ edges in $E^\phi_{m}$, $2\le \phi\le n$, also induce a path of order $p$ with pending vertices $u_{m,1}$ and $u_{m,p}$ alternating in odd and even subscripts. Note that each graph in $\mathcal L_t(p,n)$ has size $tn(p+1)$. Thus, one can check that if $p=9$, then $n=2$ so that $|\mathcal L_1(9,2)|=4$ up to isomorphism with $a^1_1=1,a^2_1=3$ or $a^1_1=1, a^2_1=2$ as in Figure~\ref{fig:4luv}. 

\begin{example} A graph in $\mathcal L_1(13,2)$ is shown in Figure~\ref{fig:l13}. The edges in $E_{1,2}$ induce an order 13 path $u_{1,1}u_{1,4}u_{1,7}u_{1,2}u_{1,5}u_{12}u_{1,9}u_{1,6}u_{1,11}u_{1,8}u_{1,3}u_{1,10}u_{1,13}$ and $F^1_1=\{u_{1,1}u_{1,5}, u_{1,13}u_{1,5}\}$ and $F^2_1=\{u_{1,1}u_{1,9}, u_{1,13}u_{1,9}\}$. \\


\begin{figure}[H]
\begin{center}
\centerline{\epsfig{file=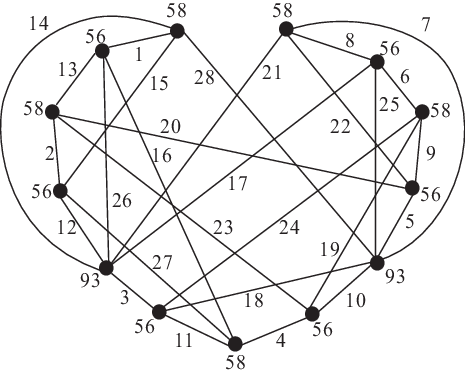, width=4cm}}
\caption{A graph in $\mathcal L_1(13,2)$ with local antimagic chromatic number 3.}\label{fig:l13}
\end{center}
\end{figure}
\end{example}

\nt Observe that each LAU-graph $G$ is tripartite so that $\chi_{la}(G)\ge 3$. In this paper, we proved that every LAU graph has local antimagic chromatic number 3. We also defined and proved that every generalized LAU graph (with components of same order and size), every circulant LAU graph (using Luv graphs of same order and size), and every generalized circulant LAU graph (with different component may have different order and size) has local antimagic chromatic number 3. Consequently, many other families of graphs  with local antimagic chromatic number 3 are also obtained.

\section{Main results}\label{sec:main}

\nt We first give an algorithmic approach to construct the $m$-th component of each graph in $\mathcal L_t(p,n)$ using $n=2$ paths, say $P_{m,\phi}$  of order $p+2$ with vertex set $V\cup \{x^\phi_{m,1},x^\phi_{m,2}\}$ for $1\le \phi\le n$. The edge set of $P_{m,1}$ is  $\{x^1_{m,1}u_{m,1}, x^1_{m,2}u_{m,p}\}\cup \{u_{m,i}u_{m,i+1}\,|\, 1\le i\le p-1\}$. For $2\le \phi\le n$, the edge set of $P_{m,\phi}$ is $\{x^\phi_{m,1}u_{m,1}, x^\phi_{m,2}u_{m,p}\}\,\,\cup$ \[\{u_{m,j^\phi_i}u_{m,j^\phi_{i+1}}\,|\, 1\le i\le p-1, |j^\phi_{i+1} - j^\phi_i| \mbox{ is odd}\ge 3,  1=j^\phi_1 \ne j^\phi_2\ne \cdots \ne j^\phi_{p-1}\ne j^\phi_p=p\}\] respectively such that every edge appears once in $P_{m,\phi}$ $(2\le \phi\le n)$.  By merging vertices $x^\phi_{m,1},x^\phi_{m,2}$ to $u_{m,2a^\phi_m+1}$ for $1\le a^1_m\ne \cdots \ne a^n_m\le \frac{p-3}{2}$, we now get the component as required. Note that the maximum value of $n$ is a function of $p$. 



\begin{theorem}\label{thm-Ln} Every LAU graph has local antimagic chromatic number 3. \end{theorem} 

\begin{proof} Recall that $\chi_{la}(G)\ge 3$ for each $G\in\mathcal L_t(p,n)$. Take $n=2$. Consider the $m$-th component $L_t(p,n|m,1,\frac{p-3}{2})$. Let the consecutive edges of $P_{m,\phi}$, $1\le \phi \le n$, be $e_{m,(\phi-1)(p+1)+i}$, $1\le i\le p+1$. Define a bijection $f : E(L_t(p,n|1,\frac{p-3}{2}))\to [1,tn(p+1)]$ such that for  $1\le i\le n(p+1)/2$, $1\le m\le t$, 

\begin{enumerate}[(1)]
\item $f(e_{m,2i}) = (m-1)n(p+1)/2+ i$;
\item $f(e_{m,2i-1}) = (2t-m+1)n(p+1)/2+1-i$. 
\end{enumerate}

\nt Observe that
\begin{enumerate}[(a)]
  \item $f^+(x_{m,1}) + f^+(x_{m,2}) =  (2t-m+1)n(p+1)/2 + (m-1)n(p+1)/2 + (p+1)/2 = (2tn+1)(p+1)/2$;
  \item $f^+(y_{m,1}) + f^+(y_{m,2}) =  (2t-m+1)n(p+1)/2 - (p+1)/2 + (m-1)n(p+1)/2 + p+1 = (2tn+1)(p+1)/2$;
  \item in $P_{m,\phi}, 1\le \phi\le n$, $f^+(u_{m,2i-1}) = (m-1)n(p+1)/2 + i + (2t-m+1)n(p+1)/2+1-i = tn(p+1)+1$ for $1\le i\le (p+1)/2$; whereas $f^+(u_{m,2i}) = (m-1)n(p+1)/2 + i + (2t-m+1)n(p+1)/2+1-(i+1) = tn(p+1)$ for $1\le i\le (p-1)/2$
'
\end{enumerate} 

\nt Combining the aboves, we can immediately conclude that in $L_t(p,n|m,1,\frac{p-3}{2})$, $f^+(u_{m,2i-1})  = tn^2(p+1)+n$ for $i\ne 2, (p-1)/2$; $f^+(u_{m,2i}) = tn^2(p+1)$ and $f^+(u_{m,3}) = f^+(u_{m,p-2}) = (2tn+1)(p+1)/2 + tn^2(p+1)+n = (2tn^2+2tn+1)(p+1)/2 + n$. Thus, $L_t(p,n|m,1,\frac{p-3}{2})$ admits a local antimagic 3-coloring. Since $L_t(p,n|m,1,\frac{p-3}{2})$ is tripartite, $\chi_{la}(L_t(p,n|m,1,\frac{p-3}{2}))=3$. 

\ms\nt In general, for $n\ge 3$, each component of every graph in $\mathcal L_t(p,n)$ can be obtained by defining a suitable $P_{m,\phi}$ $(3\le \phi\le n)$. Labeling the edges according to the function defined above, we immediately get a required local antimagic 3-coloring for every graph in $\mathcal L_t(p,n)$. This completes the proof.  \end{proof}


\begin{example}\label{eg-L2(11)} For $p=11, t=2$, using the path $P_{1,2}$ with edge set $\{x^2_{m,1}u_{m,1}, u_{m,i}u_{m,i+3}, u_{m,p}x^2_{m,2}\,|\,1\le i\le p-1\}$ (modulo $p-1$ when $i=p-2,p-1$),  we can get $\mathcal L_2(11,2) = \{L_2(11,2|1,a,b) + L_2(11,2|2,a', b')\}$ with $(a,b), (a',b')\in\{(1,2), (1,3), (1,4), (2,3), (2,4), (3,4)\}\}$. The graph $L_2(11,2|1,1,4)+L_2(11,2|2,1,3)$ is given in Figure~\ref{fig:L1113}.

\begin{figure}[H]
\begin{center}
\centerline{\epsfig{file=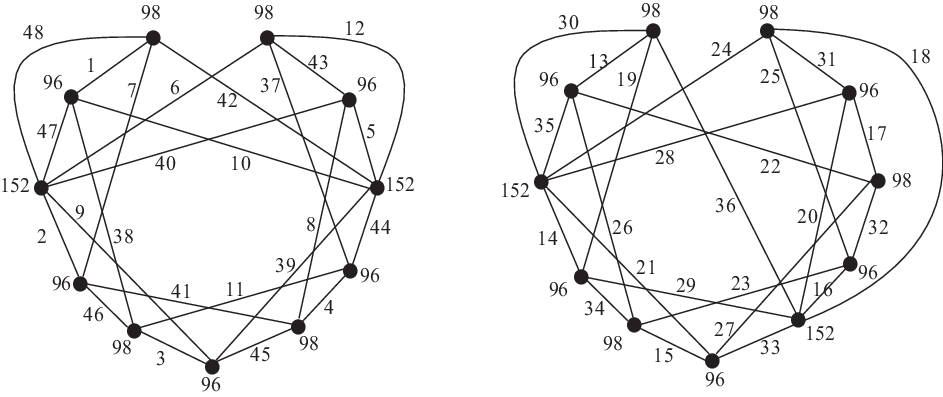, width=10cm}}
  \caption{Graph $L_2(11,2|1,1,4)+L_2(11,2|2,1,3)$ in $\mathcal L_2(11,2)$.}   \label{fig:L1113}
\end{center}
\end{figure}
\end{example}

\begin{example} Take $t = 2, p=15, n=3$, we can get various non-isomorphic 2-component graphs. One of them with the paths $P_{m,\phi}$ and the corresponding edge labels are given in Figures~\ref{fig:GL21531} and~\ref{fig:GL21532} respectively. The end vertices of each path are merged to one of the vertices with induced vertex label 97 bijectively, for example, the 3rd, 5th and 7th vertex (from the left) of degree 2 respectively. Each component of the resulting graph has two degree 8 vertices with induced vertex label 395, seven vertices of degree 6 with induced vertex label 291, and six vertices of degree 6 with induced vertex label 288.

\begin{figure}[H]
\begin{center}
\hspace{-4cm}\centerline{\epsfig{file=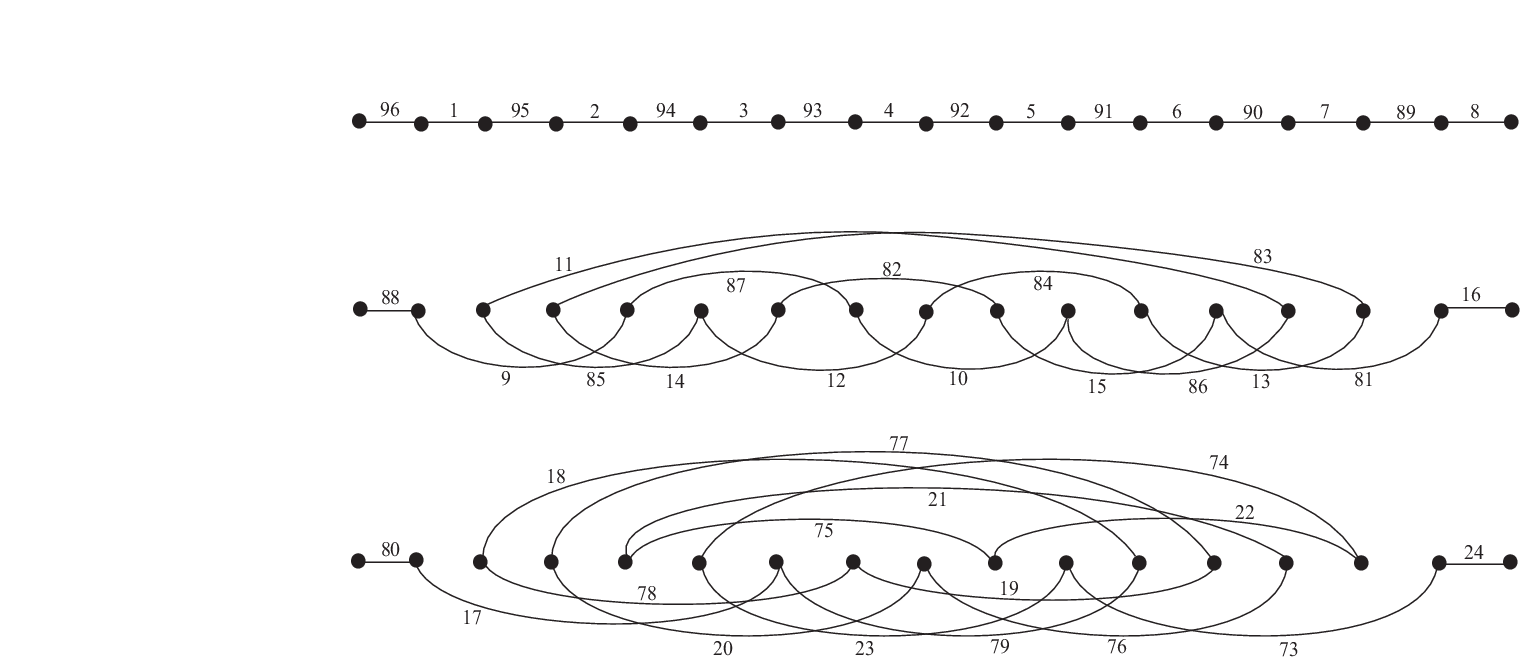, width=16cm}}
  \caption{Paths $P_{1,1}, P_{1,2}, P_{1,3}$ that give the first component.}   \label{fig:GL21531}
\end{center}
\vskip-1cm
\end{figure}
\begin{figure}[H]
\begin{center}
\hspace{-4cm}\centerline{\epsfig{file=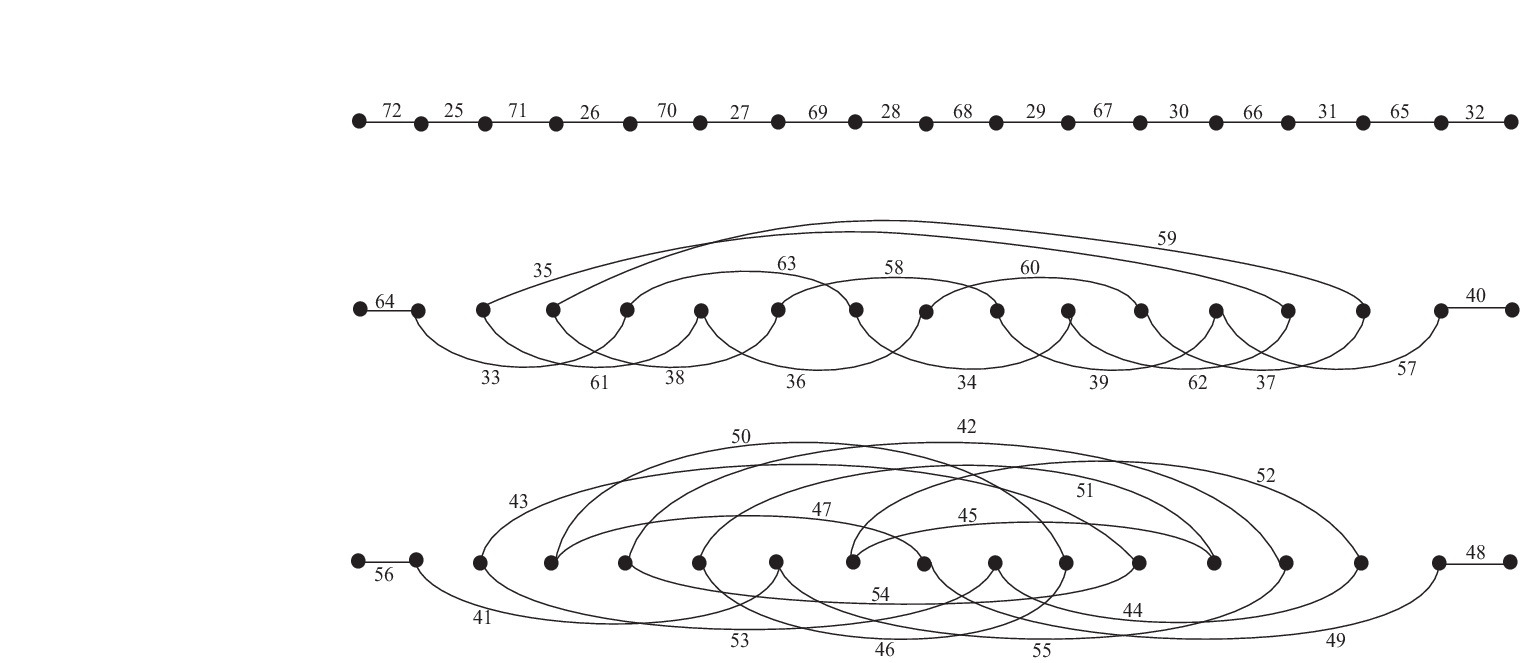, width=16cm}}
  \caption{Paths $P_{2,1}, P_{2,2}, P_{2,3}$ that give the second component.}   \label{fig:GL21532}
\end{center}
\end{figure}
\end{example} 

\nt For $9\le p_1\le p_2\le \cdots \le p_t$, each $p_i$ is odd for $1\le i\le t$, $t\ge 2$, we can now define a family of {\it generalized LAU} graph, denoted $\mathcal{GL}_t((p_1,n_1),(p_2,n_2),\ldots,(p_t,n_t))$, $n_i\ge 2, 1\le i\le t$ (also denoted $\mathcal{GL}_t(p,n)$ if $p_i=p,n_i=n$) which is a family of disconnected graphs with $t$ compoments such that the $m$-th component is $L_t(p,n_m|m,a^1_m,a^2_m,\ldots,a^{n_m}_m)$  that has $n_m$ vertices of degree $2n_m+2$ and $p_m-n_m$ vertices of degree $2n_m$. Note that $\mathcal{GL}_t(p,n) = \mathcal L_t(p,n)$.


\ms\nt {\bf Observation.} If $j_{i+1}-j_i$ is an odd constant $k\ge 3$ for all $i$, we may apply the concept of $k$-step Hamiltonian tour of cycles for odd $k$ (see~\cite[Theorem 2.5]{LLKS}) that can give us an algorithmic approach to obtain the edges in $E^\phi_{m}, 2\le \phi\le n_m$. 

\begin{example} Suppose $p_i=23$ for $1\le i\le t$, using the concept of $k$-step Hamiltonian tour, we can have a graph in $\mathcal L_t(23,n_i)$ for $k\in \{1\}\cup N$ such that $N\subseteq \{3,5,7,9\}$ using $k\in\{1,3,5,7,9\}$. Note that $k\le \frac{p_i-3}{2}$ and $2\le n_i\le |N|+1$. The corresponding paths of the $m$-th component can be as follows.
\begin{enumerate}[(i)]
\item For $k=1$, the induced path is $u_{m,1}u_{m,2}\ldots u_{m,23}$.
\item For $k=3$, the induced path is $u_{m,1}u_{m,4}u_{m,7}u_{m,10}u_{m,13}u_{m,16}u_{m,19}u_{m,22}u_{m,3}u_{m,6}u_{m,9}u_{m,12}$ $u_{m,15}u_{m,18}u_{m,21}u_{m,2}u_{m,5}u_{m,8}u_{m,11}u_{m,14}u_{m,17}u_{m,20}u_{m,23}$.
\item For $k=5$, the induced path is $u_{m,1}u_{m,6}u_{m,11}u_{m,16}u_{m,21}u_{m,4}u_{m,9}u_{m,14}u_{m,19}u_{m,2}u_{m,7}u_{m,12}$ $u_{m,17}u_{m,22}u_{m,5}u_{m,10}u_{m,15}u_{m,20}u_{m,3}u_{m,8}u_{m,13}u_{m,18}u_{m,23}$.
\item For $k=7$, the induced path is $u_{m,1}u_{m,8}u_{m,15}u_{m,22}u_{m,7}u_{m,14}u_{m,21}u_{m,6}u_{m,13}u_{m,20}u_{m,5}u_{m,12}$ $u_{m,19}u_{m,4}u_{m,11}u_{m,18}u_{m,3}u_{m,10}u_{m,17}u_{m,2}u_{m,9}u_{m,16}u_{m,23}$.
\item For $k=9$, the induced path is $u_{m,1}u_{m,10}u_{m,19}u_{m,6}u_{m,15}u_{m,2}u_{m,11}u_{m,20}u_{m,7}u_{m,16}u_{m,3}u_{m,12}$ $u_{m,21}u_{m,8}u_{m,17}u_{m,4}u_{m,13}u_{m,22}u_{m,9}u_{m,18}u_{m,5}u_{m,14}u_{m,23}$.
\end{enumerate}
Each of the paths above will have both end-vertices joining to a single vertex in $\{u_{m,i}\,|\,i=3,5,7,\ldots,21\}$ bijectively.
\end{example}

\begin{problem} Study $\chi_{la}(G)$ for $G\in \mathcal{GL}_t((p_1,n_1),(p_2,n_2),\ldots,(p_t,n_t))\setminus \mathcal{GL}_t(p,n)$. \end{problem}


\nt Observe that for graphs in $\mathcal{GL}_t(p,n)$, $t\ge 2, n\ge 2, p\ge 9$ is odd, and the defined local antimagic 3-coloring $f$, we may say every vertex of degree $2n$ is incident to $2n'$ edges (for $1\le n' < n$) with equal labels sum in $\{n'[tn(p+1)+1], n'[tn(p+1)]\}$ (respectively, of degree $2n+2$ is incident to $2n'$ edges (for $1\le n' \le n)$ with equal label sum $n'[(2tn^2+2tn+1)(p+1)/2+n]$).  We can now construct new families of connected graphs using {\it edge-swap} process as follows.

\begin{enumerate}[(1)]
\item Choose 2 vertices, say $x$ and $y$, of same induced vertex label from two different components. 
\item For each vertex, identify any $2n'$ $(1\le n' < n)$ incident edges with equal labels sum in $\{n'[tn(p+1)+1], n'[tn(p+1)]\}$ (respectively, identify any $2(n'+n'')$ incident edges with equal labels sum  $n'[tn(p+1)+1]+n''(2tn+1)(p+1)/2$ for $n'+n''\le n, 1\le n'\le n, n'' \in\{0,1\})$.
\item Redraw the $4n'$ edges (respectively, the $4(n'+n'')$ edges) so that the $2n'$ edges (respectively, the $2(n'+n'')$ edges) that are originally incident to $x$ of a component are now incident to $y$, and vice versa.
\item Repeat the above process for as many times as possible
\end{enumerate}
Consequently, by keeping all the edge labels, the new graph obtained is a connected graph that preserved the induced vertex labels. Note that all the graphs such obtained are still tripartite. Let the family of graphs such obtained using $t\ge 2$ components graph in $\mathcal{GL}_t(p,n)$ be known as {\it Circulant LAU} graphs, denoted $\mathcal{CL}_t(p,n)$. We immediately have the following theorem with the proof omitted.

\begin{theorem} Every circulant LAU graph has local antimagic chromatic number 3. \end{theorem}

\begin{example}\label{eg:CL292} Using the component $L_2(11,2|1,1,4)$ in Figure~\ref{fig:L1113}, we can get various non-isomorphic circulant LAU graphs. Two of them are given in Figures~\ref{fig:CL21121} and~\ref{fig:CL21122}.

\begin{figure}[H]
\begin{center}
\centerline{\epsfig{file=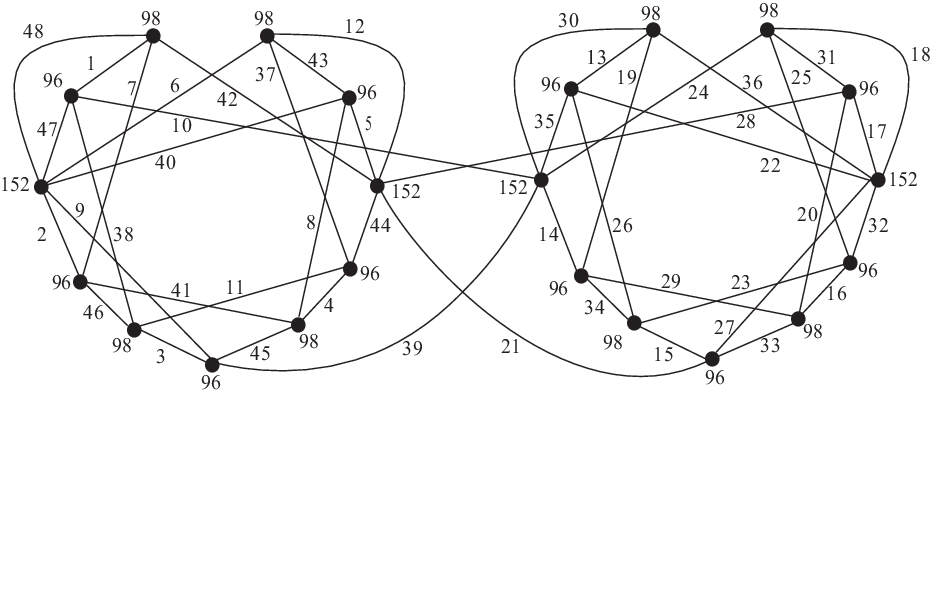, width=8cm}}
\vskip-1.8cm\caption{A graph in $\mathcal{CL}_2(11,2)$ by applying one time edge-swap.}\label{fig:CL21121}
\end{center}
\end{figure}

\begin{figure}[H]
\begin{center}
\centerline{\epsfig{file=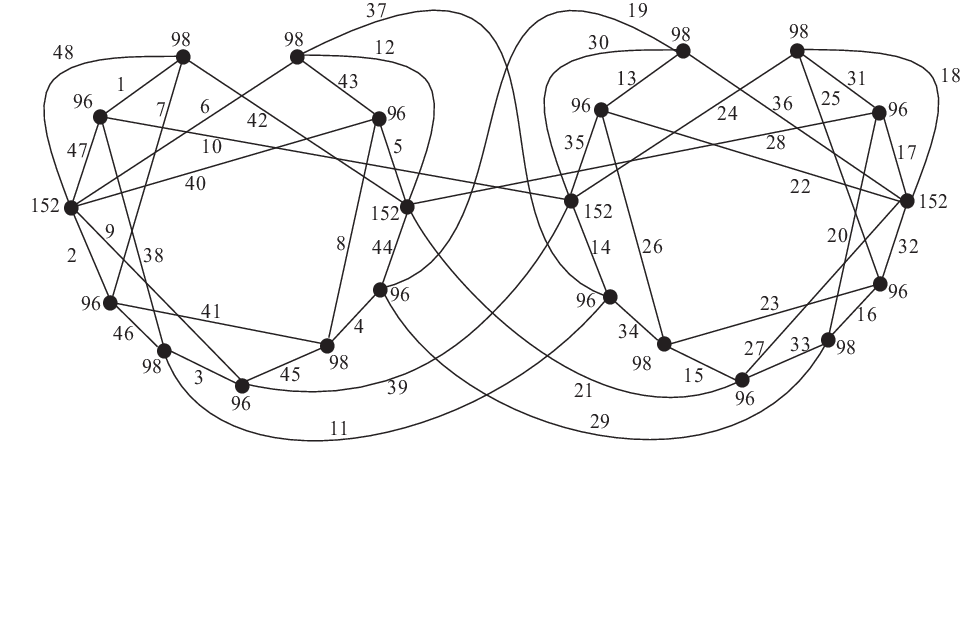, width=8cm}}
\vskip-1.6cm\caption{A graph in $\mathcal{CL}_2(11,2)$ by applying two times edge-swap.}\label{fig:CL21122}
\end{center}
\end{figure}
\end{example}

\nt More generally, we can also apply the idea of labeling graphs in $\mathcal{GL}_t(p,n)$ to define a familiy of {\it generalized circulant LAU} graphs, denoted $\mathcal{GCL}_{s_1,\ldots,s_t}(p,n)$ with $s_i\ge 2$ for at least an $i\in[1,t]$ and the $i$-th component is in $\mathcal{CL}_{s_i}(p,n)$, if $s_i\ge 2$, such that the graphs obtained also has local antimagic chromatic number 3. The theorem is stated below without proof. 

\begin{theorem} Every graph $G\in \mathcal{GCL}_{s_1,\ldots,s_t}(p,n)$ has $\chi_{la}(G)=3$. \end{theorem}

\begin{example} We give a 3-component example with $s_1=s_2=1,s_3=2$ using the two graphs in Example~\ref{eg:CL292}.

\begin{figure}[H]
\begin{center}
\centerline{\epsfig{file=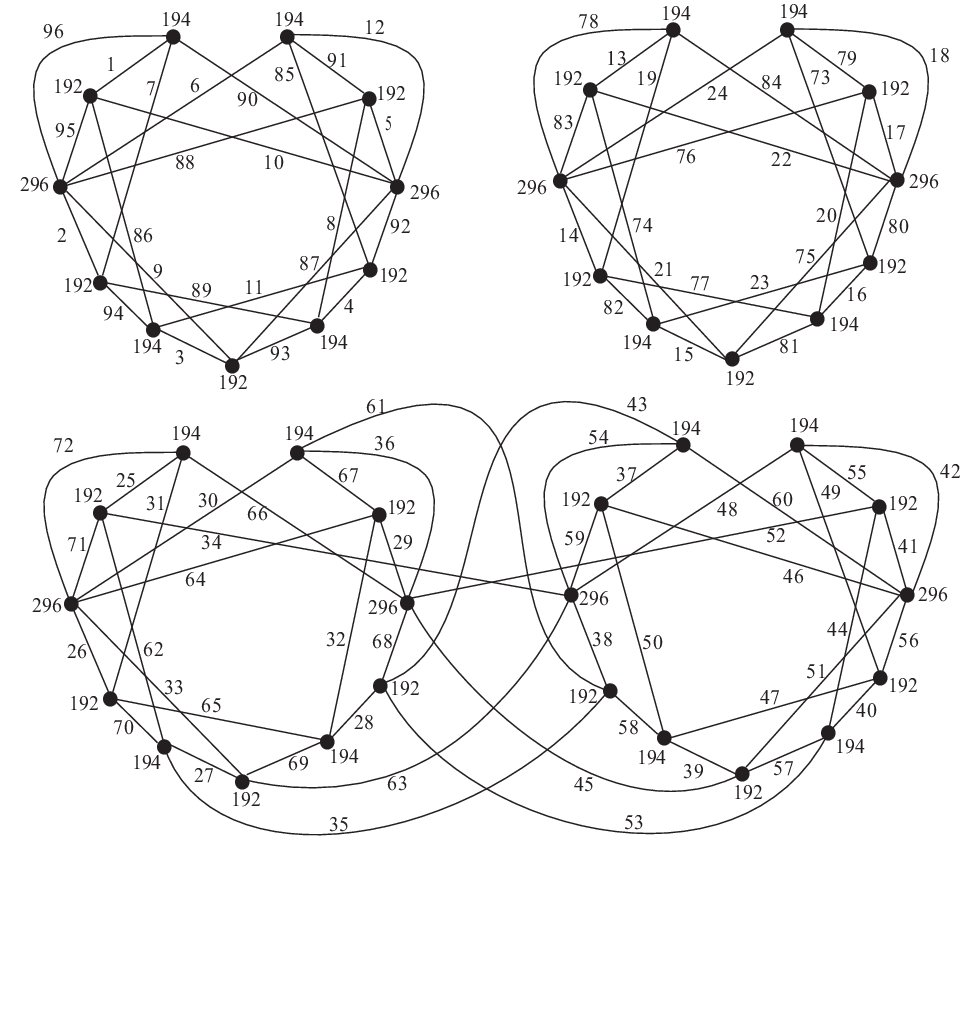, width=8cm}}
\vskip-1.8cm\caption{A 3-component graph in $\mathcal{GCL}_{1,1,2}(11,2)$.}\label{fig:GCL2112}
\end{center}
\end{figure}

\end{example}

\nt Further observe that the local antimagic 3-coloring $f$ defined in the proof of Theorem~\ref{thm-Ln} induces a 3-independent partitions of  size $tn$, $t(p-2n+1)/2$ and $t(p-1)/2$ with vertices of degree $2n+2, 2n, 2n$ respectively. Let $t=rs\ge 2$, $r\ge 1,s\ge 2$. Partition the $t$ components of each graph in $\mathcal{GL}_t(p,n) = \mathcal L_t(p,n)$ into $r$ mutually disjoint set(s), say $D_1, D_2, \ldots, D_r$, of $s$ components. For each graph in $\mathcal{GL}_t(p,n)$, we now construct three families of graphs as follows. Consider $1\le \rho\le r$.

\begin{enumerate}[(a)]
  \item Let $\mathcal{GL}^1_r(s(p-n)+n,n)$ be the family of graph(s) with $r$ component(s) such that the $\rho$-th component is of order $s(p-n)+n$ obtained by merging the $k$-th ($k\in[1,n]$) degree $2n+2$ vertex (with induced vertex label $(2tn^2+2tn+1)(p+1)/2 + n$) of each component in $D_\rho$ bijectively. Each component of the graph obtained has $n$ vertices of degree $s(2n+2)$, $s(p-2n+1)/2$ and another $s(p-1)/2$ vertices of degree $2n$.
  \item Let $\mathcal{GL}^2_r(s(p+2n-1)/2+(p-2n+1)/2,n)$ be the family of graph(s) with $r$ component(s) such that the $\rho$-th component is of order $s(p+2n-1)/2 + (p-2n+1)/2$ obtained by merging the $k$-th ($k\in [1,(p-2n+1)/2]$) degree $2n$ vertex (with induced vertex label $tn^2(p+1)+n$) of each component in $D_\rho$ bijectively. Each component of the graph obtained has $(p-2n+1)/2$ vertices of degree $2sn$, $sn$ vertices of degree $2n+2$ and another $s(p-1)/2$ vertices of degree $2n$. 
  \item Let $\mathcal{GL}^3_r(s(p+1)/2+(p-1)/2,n)$ be the family of graph(s) with $r$ component(s) such that the $\rho$-th component is of order $s(p+1)/2 + (p-1)/2$ obtained by merging the $k$-th ($k\in[1, (p-1)/2]$) degree $2n$ vertex (with induced vertex label $tn^2(p+1)$) of each component in $D_\rho$ bijectively. Each component has $(p-1)/2$ vertices of degree $2sn$, $sn$ vertices of degree $2n+2$ and another $s(p-2n+1)/2$ vertices of degree $2n$.
\end{enumerate}  

\begin{theorem}\label{thm-mergeLn} If $G\in\mathcal{GL}^1_r(s(p-n)+n,n)\cup\mathcal{GL}^2_r(s(p+2n-1)/2+(p-2n+1)/2,n)\cup\mathcal{GL}^3_r(s(p+1)/2+(p-1)/2,n)$, then $\chi_{la}(G)=3$. \end{theorem} 

\begin{proof} Clearly, $\chi_{la}(G)\ge\chi(G)=3$. Keeping the local antimagic 3-coloring for the corresponding graph as defined in the proof of Theorem~\ref{thm-Ln}. Suppose $G\in \mathcal{GL}^1_r(s(p-n)+n,n)$. Thus, $G$ admits an edge labeling such that in each component, there are $n$ degree $s(2n+2)$ vertices with induced label $s((2tn^2+2tn+1)(p+1)/2 + n)$ and the remaining degree $2n$ vertices still have induced vertex labels $tn^2(p+1)+n$ or $tn^2(p+1)$. 

\ms\nt Suppose $G\in \mathcal{GL}^2_r(s(p+2n-1)/2+(p-2n+1)/2,n)$. Thus, $G$ admits an edge labeling such that in each component, there are $(p-2n+1)/2$ degree $2sn$ vertices with induced label $s(tn^2(p+1)+n)$, the $sn$ degree $2n+2$ vertices and the remaining $s(p-1)/2$ degree $2n$ vertices still have induced vertex labels $(2tn^2+2tn+1)(p+1)/2 + n$ and $tn^2(p+1)$ respectively. 

\ms\nt Suppose $G\in \mathcal{GL}^3_r(s(p+1)/2+(p-1)/2,n)$. Thus, $G$ admits an edge labeling such that in each component, there are $(p-1)/2$ degree $2ns$ vertices with induced label $s(tn^2(p+1))$, the $sn$ degree $2n+2$ vertices and the remaining $s(p-2n+1)/2$ degree $2n$ vertices still have induced vertex labels $(2tn^2+2tn+1)(p+1)/2 + n$ and $tn^2(p+1)+n$ respectively. 

\ms\nt In each case, it is easy to check that all the induced vertex labels are distinct so that $G$ admits a local antimagic 3-coloring. Thus, $\chi_{la}(G)\le 3$. The proof is complete.
\end{proof} 

\nt Note that when $p$ is sufficiently large, it is possible that every component of a graph in $\mathcal{GL}_t(p,n)$ has the $n$ degree $2n+2$ vertices, say $u_1,\cdots,u_n$, that are not adjacent to nor having common neighbors with $n$ degree $2n$ vertices of same induced vertex label, say $v_1,\cdots,v_n$, under the local antimagic 3-coloring $f$ as defined. Let $G$ be the graph, necessarily without multiple edges nor loops, obtained by merging  $u_1,\cdots,u_n$ with $v_1,\cdots,v_n$ bijectively. We now define three more families of graphs naturally as follows.

\begin{enumerate}[(a)]
\item Let $\mathcal{GL}^4_t(p-n,n)$ be the family of graphs with $t\ge 1$ component(s) such that the $m$-th component is of order $p-n$ obtained by merging the $n$ degree $2n+2$ vertices bijectively with $n$ degree $2n$ vertices that have induced label $tn^2(p+1)+n$.
\item Let $\mathcal{GL}^5_t(p-n,n)$ be the family of graphs with $t\ge 1$ component(s) such that the $m$-th component is of order $p-n$ obtained by merging the $n$ degree $2n+2$ vertices bijectively with $n$ degree $2n$ vertices that have induced label $tn^2(p+1)$.
\item Let $\mathcal{GL}^6_t(p-n+1,n)$ be the family of graphs with $t\ge 1$ component(s) such that the $m$-th component is of order $p-n+1$ obtained by merging all the $n$ degree $2n+2$ vertices into a vertex of degree $n(2n+2)$.
\end{enumerate}

\begin{theorem}\label{thm-mergeLn2} If $G\in \mathcal{GL}^4_t(p-n,n)\cup \mathcal{GL}^5_t(p-n,n) \cup \mathcal{GL}^6_t(p-n+1,n)$, then $\chi_{la}(G)=3$. \end{theorem}

\begin{proof} By definition, $G$ is still tripartite so that $\chi_{la}(G)\ge \chi(G) = 3$. (a) Note that $u_i$ $(1\le i\le n)$ has induced vertex label $(2tn^2+2tn+1)(p+1)/2+n$. If $v_i$ $(1\le i\le n)$ has induced vertex label $tn^2(p+1)+n$, after merging, the degree $4n+2$ vertices obtained has induced vertex labels $(4tn^2+2tn+1)(p+1)/2+2n$ which is larger than the remaining degree $2n$ vertices with induced vertex labels $tn^2(p+1)+n$ and $tn^2(p+1)$. (b) Similarly if  $v_i$ has induced vertex label $tn^2(p+1)$. (c) The degree $n(2n+2)$ vertex has induced vertex label $n(2tn^2+2tn+1)(p+1)/2+n^2$ while the other vertices have the same induced vertex labels. Thus, $\chi_{la}(G)\le 3$. This completes the proof. \end{proof}

\begin{example} A graph in $\mathcal L_1(17,2)$ that can give a graph in Theorem~\ref{thm-mergeLn2} is shown in Figure~\ref{fig:L17}.
\begin{figure}[H]
\begin{center}
\hspace{2cm}\epsfig{file=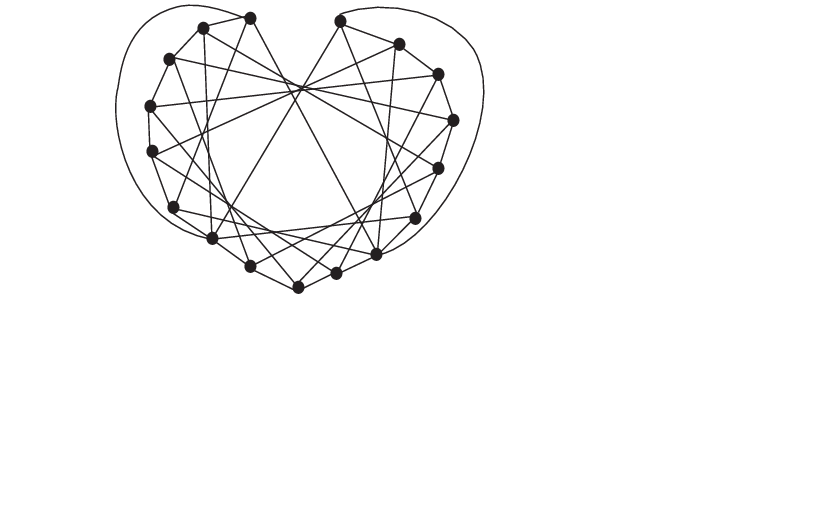, width=8cm}
\vspace{-2cm}
  \caption{A graph in $\mathcal L_1(17,2)$.}\label{fig:L17}
\end{center}
\end{figure}
\end{example}

\begin{remark}  We note that the construction of graphs in Theorems~\ref{thm-mergeLn} and~\ref{thm-mergeLn2} using the graphs in $\mathcal{GL}_t(p,n)$ can be done by using the graphs in $\mathcal{CL}_t(p,n) \cup \mathcal{GCL}_{s_1,\ldots,s_t}(p,n)$ as well so that theorems similar to Theorems~\ref{thm-mergeLn} and~\ref{thm-mergeLn2}  can be obtained too.  \end{remark}

\section{Conclusion and open problems} 

\nt In this paper, we introduced algorithmic approach to contruct various families of non-regular (possibly disconnected) tripartite graphs of even size and proved that all these graphs have local antimagic chromatic number 3. 

\ms\nt In~\cite[Lemmas 2.3 \& 2.4]{LSN}, the authors obtained sufficient conditions to have $\chi_{la}(G)=\chi_{la}(G-e)$. Let $e$ be an edge of a graph $G$ in any of the theorems in Section~\ref{sec:main} with label $1$ or $|E(G)|$ under the local antimagic 3-coloring $f$ as defined. One may check the conditions in~\cite[Lemmas 2.3 \& 2.4]{LSN} to obtain the exact value of $\chi_{la}(G-e)$. In general, we can have the following problem.






\begin{problem} If $e$ is an edge of $G$, determine $\chi_{la}(G-e)$. \end{problem}

\nt Observe that if we redefine graphs in $\mathcal{GL}_t(p,n)$ so that the degree $2n+2$ vertices are $u_{m,i}, u_{m,i'}$ where $i,i'$ not both odd for  $1\le m\le t$, then we have either a bipartite graph $\in\mathcal B$ with distinct partite set size or else a tripartite graph $T\in\mathcal T$ that admits a local antimagic 3- or 4-coloring respectively. Thus, $2\le \chi_{la}(B)\le 3$ and $3\le \chi_{la}(T)\le 4$.

\begin{problem} Determine $\chi_{la}( B)$ and $\chi_{la}(T)$. \end{problem} 

\nt We end this paper with the following problems.

\begin{problem} Determine the maximum value of $n$ for each possible value of $p$.  \end{problem} 

\begin{problem} Study the properties of the various graph polynomials (and uniqueness), and graph parameters of every graph in Section~\ref{sec:main} such as (but not limited to) chromatic, domination, independent, star (or neighborhood) polynomials, magicness, antimagicness, Roman domination number and Sudoku number.  \end{problem}

\section{Statements and Declarations}

\noindent No funds, grants, or other support was received. The author has no relevant financial or non-financial interests to disclose.

\end{document}